# An Inverse Gaussian Process Monte Carlo algorithm for estimation and uncertainty assessment of hydrologic model parameters

Jiangjiang Zhang[1], Weixuan Li[2]

[1] College of Environmental and Resource Sciences, Zhejiang University, Hangzhou, 310058, China

[2] Pacific Northwest National Laboratory, Richland, WA 99352, USA

**Abstract.** Solving hydrologic inverse problems usually requires repetitive forward simulations. One approach to mitigate the computational cost is to build a surrogate model, i.e., an approximate mapping from model parameters (input) to observable quantities (output), so the forward simulations can be done quickly. Alternatively, if the surrogate is constructed to approximate the inverse mapping from model outputs to parameters, the parameter estimates can be obtained directly by treating measurements as inputs to this inverse surrogate. Moreover, the uncertainties of parameters can be quantified by propagating the measurement uncertainties in a straightforward Monte Carlo manner. Based on this idea, we proposed a novel surrogate-based approach for parameter estimation and uncertainty assessment, i.e., the Inverse Gaussian Process Monte Carlo (IGPMC) algorithm. The Gaussian Process (GP) regression is used to directly approximate the inverse function of the model output-input relationship. For ill-posed problems, i.e., when there exist non-unique sets of input parameters that all produce an identical system output, multiple inverse GP systems are constructed and multiple parameter estimates can be obtained accordingly. The accuracy and efficiency of this IGPMC algorithm were demonstrated through four numerical case studies. Results obtained from the Markov Chain Monte Carlo (MCMC) are used as references to assess our new proposed method. It was shown that, the IGPMC algorithm can generally obtain reliable parameter estimates with an affordable computational cost.

# 1. Introduction

Many parameters in hydrological models are conceptual, which are difficult or sometimes impossible to be directly measured. Spatial heterogeneity of the model parameters further increases the difficulty in their measurement. Thus accurate parameter estimation is critical for numerical modeling in hydrologic sciences.

On the other hand, many model state variables such as stream flow for rainfall-runoff models, hydraulic head and concentration for groundwater flow and solute transport models, can be measured directly or even monitored continuously in situ. With measurements of the state variables, model parameters can be estimated indirectly by solving an inverse problem [*Carrera et al.*, 2005; *Matott et al.*, 2009; *Oliver and Chen*, 2011; *Tartakovsky*, 2013; *Zhou et al.*, 2014].

The deterministic inverse methods, such as gradient-based methods [*Gupta and Sorooshian*, 1985], and genetic algorithm [*El Harrouni et al.*, 1996], are usually accomplished by solving least-squares minimization problems by searching for a single parameter set that best fits the measurements. However, there might be multiple local optima in the parameter space and the ability to find the global optimum is essential for these methods [*Pan and Wu*, 1998]. The key problem for the deterministic methods is that, in these approaches, the system and parameter uncertainties cannot be fully characterized [*Ma and Zabaras*, 2009].

To take the system and parameter uncertainties into account, it is more reasonable to formulate a stochastic description of the unknown parameters with given measurement data. Among the stochastic inverse methods, the Bayesian method is one of the most popular approaches to estimate the unknown parameters. Except for special cases where analytical forms can be obtained [*Woodbury and Ulrych*, 2000], posterior distributions of parameters are usually estimated with sampling methods. In the Bayesian framework, the information gained from measurements can be well assessed in the posterior distributions of parameters.

The Markov Chain Monte Carlo (MCMC) method can be used as an ideal parameter estimation method for highly nonlinear and non-Gaussian problems involving complex processes [*Andrieu et al.*, 2003; *Vrugt et al.*, 2003]. For its general applicability, MCMC method is receiving increasing popularity in hydrologic sciences [*Laloy et al.*, 2013; *Liu et al.*, 2010; *Smith and Marshall*, 2008; *Vrugt et al.*, 2009a; *Vrugt et al.*, 2008; *Vrugt et al.*, 2009b; *Zeng et al.*, 2012; *Zhang et al.*, 2015; *Zhou et al.*, 2014]. However, a large number of model evaluations are usually needed to sufficiently explore the posterior parameter space. Even with some advanced MCMC algorithms such as the Delayed Rejection and Adaptive Metropolis (DRAM) algorithm [*Haario et al.*, 2006] and the Differential Evolution Adaptive Metropolis (DREAM) algorithm [*Vrugt et al.*, 2008; *Vrugt et al.*, 2009b], tens of thousands of model evaluations are usually required.

One possible approach to accelerate MCMC is to construct a surrogate for the original model through polynomial chaos expansion [*Laloy et al.*, 2013] or sparse grid interpolation [*Zeng et al.*, 2012; *Zhang et al.*, 2013; *Zhang et al.*, 2015]. Then this surrogate is used to replace the original system in the likelihood evaluation. However, when the number of parameters is large and the model is highly nonlinear, it is unrealistic to construct an accurate surrogate.

Gaussian Process (GP) is a stochastic process assuming that the realizations have Gaussian distributions, i.e., the process is specified by its mean and covariance functions. GP has been widely used in many aspects associated with machine learning such as regression and classification [*Rasmussen and Nickisch*, 2010]; applied in geostatistics for kriging interpolation to estimate the spatial distribution of geological property; and used to approximate the model input-output relationship. The latter was achieved by fitting this stochastic process with certain model inputs and corresponding outputs [*Kennedy and O'hagan*, 2001; *Razavi et al.*, 2012; *Sun et al.*, 2014; *Li, 2014*].

In contrast to previous applications, this research employed GP to construct surrogate systems mapping from model outputs to model parameters. In other words, GP was used to approximate the inverse function of the system model to estimate the model parameters by using the measurements as inputs. Furthermore, the uncertainty of the measurements that propagated to model parameters through this inverse GP system was assessed in a straightforward Monte Carlo (MC) manner.

The paper is organized as follows: The methods are formulated in Section 2. In Section 3, the performance of proposed methods is illustrated with four synthetic numerical case studies. Finally, some conclusions are provided in Section 4.

## 2. Methods

In a hydrologic model, measurements $\mathbf{d}$ can be expressed as

$$\mathbf{d} = F(\mathbf{m}) + \boldsymbol{\varepsilon}, \tag{1}$$

where $\mathbf{m}$ and $F(\mathbf{m})$ are $n_m \times 1$ and $n_d \times 1$ vectors of the model parameters and outputs, respectively, and $n_m$ and $n_d$ are the dimensions of parameters and measurements, respectively; $\boldsymbol{\varepsilon}$ is a $n_d \times 1$ vector of measurement errors with certain distribution. We are interested in estimation and uncertainty assessment of the model parameters $\mathbf{m}$ from noisy measurements $\mathbf{d}$. In this paper, we propose a simple while effective method, i.e., the Inverse Gaussian Process Monte Carlo algorithm, to estimate the model parameters.

### 2.1. Gaussian Process Regression

The Inverse Gaussian Process Monte Carlo algorithm is based on the construction of a system mapping from model outputs to parameters with GP regression, which is different from the traditional use of GP to construct a surrogate system mapping from model parameters to outputs.

The following is a brief introduction for the traditional use of GP regression to approximate the mapping from model parameters $\mathbf{m}$ to original model outputs $F(\mathbf{m})$. The fundamental philosophy of GP is that each element of $F(\mathbf{m})$ is assumed to be a realization of Gaussian stochastic process $G(\mathbf{m})$, which can be specified by its first two order statistics,

$$G(\mathbf{m}) \sim N(\boldsymbol{\mu}(.), C(.,.)), \tag{2}$$

where $\mu(.)$ and $C(.,.)$ are the mean and covariance functions for the GP process.

Since we cannot evaluate $F(\mathbf{m})$ at every point of parameter $\mathbf{m}$, the outputs of GP are associated with uncertainty. Suppose we already have outputs for $N$ sets of parameter base points $\{\mathbf{m}_B^{(1)}, \mathbf{m}_B^{(2)}, ..., \mathbf{m}_B^{(N)}\}$, then the system outputs can be conditioned to these base points. This conditioned process, demoted as $G_{|\mathbf{B}}(\mathbf{m})$, is still a GP with its mean and variance at an arbitrary parameter point $\mathbf{m}$, given as

$$\mu_{|\mathbf{B}}(\mathbf{m}) = \mu(\mathbf{m}) + C_{\mathbf{mB}} C_{\mathbf{BB}}^{-1} (F(\mathbf{B}) - \mu(\mathbf{B})), \tag{3}$$

and

$$\sigma^2_{|\mathbf{B}}(\mathbf{m}) = C(\mathbf{m}, \mathbf{m}) - C_{\mathbf{mB}} C_{\mathbf{BB}}^{-1} C_{\mathbf{Bm}}, \tag{4}$$

where $C_{\mathbf{mB}}$ is a $1 \times N$ vector whose $i$th element is $C(\mathbf{m}, \mathbf{m}_B^{(i)})$, $C_{\mathbf{BB}}$ is a $N \times N$ matrix with $C(\mathbf{m}_B^{(i)}, \mathbf{m}_B^{(j)})$ as its $i$th row and $j$th column element, $C_{\mathbf{Bm}}$ is the transpose of $C_{\mathbf{mB}}$. It can be seen that, the variances of the system outputs are zero at the base points, while the variance at any other parameter point is also reduced because of the correlation between $\mathbf{m}$ and base points, as shown in Eq. (4). Since it requires evaluating the original model once at each chosen base point, one key to reduce the computational cost is the wise selection of base points.

The performance of GP is also influenced by the choice of mean and covariance functions in Eqs. (3) and (4). There are many candidate mean and covariance functions [*Rasmussen and Nickisch*, 2010], which could be determined based on model outputs at the base points. One popular method proposed by [*Jones et al.*, 1998] is to find the one that maximizes a likelihood function of the model outputs at these base points. The likelihood function is with $2k$ parameters, where $k$ is the dimension of inputs. As this method resorts to solving an optimization problem with $2k$ parameters, the computational cost for determining the function forms would be very high if $k$ is large.

## 2.2. The Inverse Gaussian Process Monte Carlo Algorithm

In this section, an efficient parameter estimation method based on the construction of inverse GP system is proposed. The idea behind this algorithm is very straightforward as illustrated in the following.

In traditional approaches, with selected base points of model parameters and corresponding model outputs, the GP system can be constructed to approximate the model input-output relationship. Intuitively, if an inverse GP system is constructed from model outputs to the model parameters, we can obtain the parameter estimates directly by using measurements as inputs to the inverse GP system. In other words, the inverse GP system is essentially a surrogate for the inverse function of the original system model.

However, as is well known, the inverse function may not even exist if multiple inputs can result in an identical system output. This non-unique mapping also gives rise to the ill-posedness for the traditional optimization-based inverse methods, or multi-modal posterior distributions from a Bayesian point of view. Hence direct construction of the inverse GP system may be inappropriate. To cope with this problem, we can divide the base points into different groups with clustering analysis through analyzing the input-output realizations. Then multiple inverse GP systems can be respectively constructed based on each group of realizations. In this way, multiple parameter estimates can be obtained simultaneously from those inverse GP systems, and the ill-posedness of inverse problem can be alleviated.

The actual measurements $\mathbf{d}$ can be expressed as true model outputs with additive measurement errors, as shown in Eq. (1). Without losing generality, the measurement errors in this work are assumed to follow Gaussian distributions with zero means, i.e., $\varepsilon \sim N(\mathbf{0}, \boldsymbol{\sigma}^2)$, where $\boldsymbol{\sigma}^2$ is the error covariance matrix. Thus, the true model outputs can be expressed as $F(\mathbf{m}) = \mathbf{d} - \varepsilon$, i.e., $F(\mathbf{m}) \sim N(\mathbf{d}, \boldsymbol{\sigma}^2)$. By using measurement samples generated from $N(\mathbf{d}, \boldsymbol{\sigma}^2)$ as inputs to the inverse GP system(s), we can directly obtain samples of parameters, and from which any wanted statistics can be obtained.

As described above, this is actually a Monte Carlo method, which is widely used for uncertainty quantification. It should be noted that, the applicability of this method relies on two assumptions: the first is that the clustering analysis can successfully alleviate the possible ill-posedness of inverse problem; the second is that the GP system can approximate the inverse function, which is treated as a black box.

The schematic diagram of the IGPMC algorithm is shown in Figure 1. The detailed implementation is described as follows.

[Figure 1]

(1) The IGPMC algorithm starts with drawing initial candidate base points from prior distribution of parameters. Here $N$ random parameter samples $\mathbf{m}^{(i)}$ and corresponding model outputs $F(\mathbf{m}^{(i)})$ are obtained, where $i = 1, 2, ..., N$.

For the IGPMC algorithm, the inputs for the GP system are actually model responses. As the dimension of $F(\mathbf{m})$ is usually very high, determining the mean and covariance functions with the method proposed by [*Jones et al.*, 1998] is usually difficult. Here we propose an empirical approach which is computationally efficient to choose the mean and covariance functions. It is reasonable to choose the mean values of selected base points as the mean function $\mu(F(\mathbf{m}))$. The covariance function $C(.,.)$ has the simple form as below,

$$C(F(\mathbf{m}^{(1)}), F(\mathbf{m}^{(2)})) = \sigma_B^2 \exp[-\sum_{i=1}^{n_d} \alpha_i |F(\mathbf{m}^{(1)}) - F(\mathbf{m}^{(2)})|^2], \tag{5}$$

where $\sigma_B^2$ is the variance of model outputs at selected base points, $\alpha_i$ is the correlation length parameters.

The inverse of correlation length parameter $\frac{1}{\alpha_i}$ controls the distance beyond which two points have relatively weak correlation in the *i*th input dimension. For simplicity and computational efficiency, in this paper, as a rule a thumb, we assume that $\frac{1}{\alpha_i} = q \cdot \sigma_i$, where $\sigma_i$ is the standard deviation of $F(\mathbf{m})_i, i = 1, 2, ..., n_d$, and $q$ is chosen so that the mean value of $\frac{C_{mB}}{\sigma_B^2}$ is around 0.75-0.95. Therefore, the choice of covariance function is greatly simplified

(2) To obtain more accurate inverse GP system efficiently, the base points used for constructing the inverse GP system are selected from the candidates by comparing the weighted distance $D(\mathbf{m}^{(i)})$ between $F(\mathbf{m}^{(i)})$ and actual measurements $\mathbf{d}$,

$$D(\mathbf{m}^{(i)}) = \frac{1}{2}(F(\mathbf{m}^{(i)}) - \mathbf{d})^T (\sigma^2)^{-1} (F(\mathbf{m}^{(i)}) - \mathbf{d}), \tag{6}$$

where $\sigma^2$ is the error covariance matrix. Smaller values of $D$ indicate better base points. The $K$ ($K \leq N$) parameter samples $\mathbf{m}^{(k)}$, $k = 1, 2, ..., K$, with smaller $D$ values are chosen as base points for the inverse GP system. Although $D(\mathbf{m}^{(i)})$ in Eq. (6) is suggested here, other measures of mismatch between model outputs at base points and measurements can be adopted. The choice may refer to formal or informal likelihood functions used in MCMC [*Schoups and Vrugt*, 2010; *Zhang et al.*, 2013].

(3) To cope with the possible ill-posedness of inverse problem, the base points selected at Step (2) are divided into several groups by clustering analysis (e.g., agglomerative hierarchical clustering analysis [*Day and Edelsbrunner*, 1984] and *K*-means clustering analysis [*Hartigan and Wong*, 1979]). The idea behind this treatment can be explained as follows. After Step (2), corresponding system outputs of chosen parameter samples are relatively close to the measurements. If the ill-posedness of inverse problem exists, i.e., parameters within several ranges correspond to a similar system output, then these parameter samples can be divided into different groups according to the distance between each other.

In our practical implementation of the IGPMC algorithm, the agglomerative hierarchical clustering analysis is used [*Day and Edelsbrunner*, 1984]. A relatively large preliminary cluster number $q$ (e.g., 5) is chosen. After the analysis, the base points are preliminarily divided into $q$ groups. Through checking the number of samples and $D$ values in each groups, the groups with typically small numbers of samples and bigger $D$ values are discarded, i.e., the rest $p$ groups of base points are kept.

(4) With these $p$ groups of parameter base points, inverse GP

systems $G_{(l)}, l = 1, 2, ..., p$, are respectively constructed from model outputs to parameters according to Eqs. (3-5). One realization of measurements ($\mathbf{d}^* = \mathbf{d} + \boldsymbol{\varepsilon}$) is used as inputs for each inverse GP system to obtain $p$ estimates of the model parameters $G_{(l)}(\mathbf{d}^*), l = 1, 2, ..., p$. Meanwhile, random perturbations $\boldsymbol{\varepsilon}$ with the same statistics as measurement errors are added to the model outputs at selected base points.

(5) If predefined stopping criteria are satisfied, the above procedure stops, and the algorithm goes to Step (6) for the assessment of parameter uncertainty. Otherwise it just adds the $p$ new sets of parameter samples and their corresponding model outputs into the pool of candidate base points (the pool size increases by $p$), then go back to Step (2). The stopping criteria can be set as the maximum allowed number of model evaluations. The convergence diagnosis can be implemented through visually checking the sequential samples generated during this process.

(6) Based on the final inverse GP system(s) constructed in Step (5), the uncertainty of parameters can be evaluated by propagating the uncertainty of measurements, as illustrated in Figure 1. This can be realized through a straightforward Monte Carlo manner, since the inverse function is already available. $M$ realizations of measurements can be generated from $N(\mathbf{d}, \boldsymbol{\sigma}^2)$. Then these realizations are treated as inputs to the inverse GP system(s) to obtain corresponding parameter samples.

Finally, the statistics of model parameters can be obtained by analyzing these parameter samples. It should be noted that, although $M$ can be very large, no further original model evaluation is needed. Therefore, the computational cost at this step is negligible.

## 3. Case Studies

We demonstrate the applicability of the Inverse Gaussian Process Monte Carlo algorithm in the following four case studies with increasing model complexity. The posterior distributions of parameters obtained by MCMC were used as references. In this paper, an efficient MCMC algorithm known as the Differential Evolution Adaptive Metropolis algorithm, DREAM$_{(ZS)}$ [*Vrugt et al.*, 2008; *Vrugt et al.*, 2009b] was used.

### 3.1. Case Study 1 : A Simple Bimodal Case

In this case, the applicability of the IGPMC algorithm in an ill-posed inverse problem was investigated. The following simple function is considered,

$$d = m^2 + \varepsilon, \qquad (7)$$

where d is a one-dimensional measurement, m is a scalar parameter, $\varepsilon \sim N(0, 0.01^2)$ is an additive Gaussian measurement error.

Let the prior $p(m)$ be a uniform distribution with range [-1, 1], the true parameter $m^* = 0.230$, and the noisy measurement $d = 0.0414$. It is easy to see that, when the parameter value is around -0.230, the model output is also close to the measurement value. Therefore, the inverse solution is not unique. In the Bayesian point of view, the posterior distribution of parameter will be bimodal.

Figure 2 presents the distributions of $m$ inferred from the posterior samples generated by DREAM$_{(ZS)}$ algorithm, and samples generated by the IGPMC algorithm. There were two parallel chains for DREAM$_{(ZS)}$, and the maximum number of total model evaluations was 2,000. The initial 500 samples were discarded and the rest were used to estimate the posterior distribution of $m$. The likelihood for DREAM$_{(ZS)}$ was Gaussian likelihood with homoscedastic measurement error. While for the IGPMC algorithm, the number of random parameter samples drawn from prior distribution was 500, and 400 iterative steps were then implemented. The initial number of clusters was chosen as 5, and then those clusters with less than 30 samples and bigger $D$ values were discarded, based on which multiple inverse GP systems were constructed simultaneously. For the IGPMC algorithm, totally 1,600 model evaluations were required.

After those 1,600 model evaluations, 1,000 perturbed measurement samples were used as inputs to the two inverse GP systems, then the resulted 2,000 samples of $m$ were used for analysis. It should be noted that, the generation of 2,000 parameter samples did not require any extra original model evaluations. It is shown in Figure 2 that, the bimodal conditional parameter distribution is well identified by both the DREAM$_{(ZS)}$ and IGPMC algorithm.

[Figure 2]

## 3.2. Case Study 2 : The HYMOD Model

This case study demonstrates the efficiency of the IGPMC algorithm in estimating the parameters of HYMOD [*Moore*, 1985], a classical five-parameter conceptual rainfall-runoff model, which was also used in the numerical studies for DREAM$_{(ZS)}$ in [*Vrugt et al.*, 2003; *Vrugt et al.*, 2008].

The HYMOD model was developed to represent the hydrological processes within a watershed, which consists of a simple rainfall excess model connected with two series of linear reservoirs [*Vrugt et al.*, 2003]. There are five parameters for HYMOD in this case study, the maximum storage capacity in the watershed, $C_{max}[L]$, the degree of spatial variability of the soil moisture capacity within the watershed, $b_{exp}$, the factor distributing the flow between the two series of reservoirs, $\alpha^*$, the residence time of the linear slow reservoirs, $R_s[T]$, and the residence time of the quick reservoirs, $R_q[T]$. The prior ranges and true values for the five parameters are listed in Table 1. The stream flow measurements used for parameter estimation was generated with the true parameters with additive measurement errors $\varepsilon \sim N(0, \sigma^2)$, where $\sigma$ is a vector that is equal to 10% the value of model outputs given the true model parameters. [*Moore*, 1985] provided the detailed description of HYMOD.

[Table 1]

In this case, to provide a reference of parameter distributions, DREAM$_{(ZS)}$ with three parallel chains and altogether 7,000 model evaluations was implemented. Gaussian likelihood with heteroscedastic measurement errors was used. The convergence was reached with about 3,000 model evaluations, and the last 4,000 samples were used to estimate the posterior distributions of parameters. For the IGPMC algorithm, the initial number of random parameter sets was 1,000, and 2,000 iterative steps were employed, i.e., the total number of model evaluations was 3,000.

[Figure 3]

Figure 3 shows the trace plots of the five parameters of HYMOD generated by the IGPMC algorithm (i.e., all the candidate base points sequentially generated during the inverse GP construction). The IGPMC algorithm converged to the true values with about 500 iterative steps after the initial 1,000 random parameter sample drawings. For the IGPMC algorithm, to assess the uncertainty of parameters, 1,000 samples of parameters were generated with direct MC sampling as described at Step (6) in section 2.2. It should be noted here that, although 1,000 samples of parameters were obtained, no original model evaluation was needed, and the time needed for this MC process was negligible.

As shown in Figure 4, compared with DREAM$_{(ZS)}$, the IGPMC algorithm can obtain comparable distributions of the HYMOD parameters. Moreover, the distributions obtained by the IGPMC algorithm (represented by red dashed lines in Figure 4) are smoother.

[Figure 4]

The uncertainty assessment of parameters in the IGPMC algorithm is based on the propagation of measurement uncertainty to parameter uncertainty through the inverse GP system. Therefore, the statistics of measurement errors are essential to the IGPMC algorithm. There are many cases that the statistics of measurement errors are unknown a priori. In DREAM$_{(ZS)}$, the unknown statistics of measurement errors can be estimated together with parameters [*Vrugt et al.*, 2008]. Similarly, we proposed a method to estimate the statistics of measurement errors if they are not known a priori.

In each iterative step, the residuals between model outputs at selected base points and measurements can be used to infer the statistics of measurement errors. For example, if the measurement errors are assumed to follow Gaussian distributions, the mean and standard deviation of these residuals can be used to represent the measurement error statistics. Thus, the parameters and statistics of measurement errors can be estimated simultaneously in IGPMC.

To illustrate the applicability of this treatment, we applied IGPMC algorithm and DREAM$_{(ZS)}$ (Gaussian likelihood with measurement integrated out was used) for this case without the information of measurement errors, and the parameter distributions versus the true values are shown in Figure 5. It is clearly shown that, both the IGPMC and DREAM$_{(ZS)}$ algorithms can provide acceptable parameter estimation results when the measurement error information is not known a priori. However, compared with the case where the statistics of measurement errors are known (Figure 4), the maximum-a-posterior (MAP) estimates of model parameters are more likely to deviate from the true values, and the variances of parameters are larger. The comparison between the actual measurement errors and estimated measurement errors (the mean values of the residuals between model outputs at selected base points and measurements) are shown in Figure 6. It is shown that, in this synthetic case study, the added measurement errors can be well estimated.

[Figure 5]

[Figure 6]

## 3.3. Case Study 3 : A Contaminant Source Identification Case

In this case study, we tested the Inverse Gaussian Process Monte Carlo algorithm for a contaminant source identification problem in steady saturated flow.

As shown in Figure 7, the flow domain is $20 [L]$ in *x* direction and $10 [L]$ in *y* direction. The upper and lower boundaries are no flow boundaries, while the left and right boundaries are constant head boundaries with pressure heads of $12 [L]$ and $11 [L]$, respectively. In this case study, steady water flow was considered, and the flow equation was solved numerically with MODFLOW [*Harbaugh et al.*, 2000]. The solute transport was calculated numerically with MT3DMS [*Zheng and Wang*, 1999]. From 1 [T] to 6 [T], a contaminant source **S** was released with a potential area denoted by a red dashed rectangle.

[Figure 7]

In this case, the porosity and the dispersivities were assumed to be known with values as porosity $\theta = 0.25$, the longitudinal dispersivity $\alpha_L = 0.3 [L^2 T^{-1}]$ and the transverse dispersivity $\alpha_T = 0.03 [L^2 T^{-1}]$, respectively. The conductivity field was simplified with three hydraulic zones. In each zone, the hydraulic conductivity $K_i [LT^{-1}]$ (represented by its log value, i.e., $Y_i = Log K_i, i = 1, 2, 3$) was assumed to be homogenous and its value to be unknown. The contaminant source was described by 8 parameters, including the location $(x_s, y_s)$ and time-varying strengths $S_{si} [MT^{-1}]$ for $t_i = i [T] : (i+1) [T], i = 1, 2, ..., 6$. The main goal of this case study was to identify the unknown contaminant source parameters and conductivity value for each zone.

As shown in Table 2, the prior distributions of the 11 unknown parameters (8 source parameters and 3 conductivity parameters) were assumed to be uniform with given ranges. To estimate these parameters, measurements for the concentration and hydraulic head at 5 locations were available (represented by blue dots in Figure 7). The concentration measurements were obtained every $2[T]$ from $t = 2[T]$ to $t = 10[T]$. Since the flow was steady, 5 hydraulic head measurements were obtained only once at these measurement locations. The measurement errors for concentration and head were assumed to follow $N(0, 0.05^2)$ and $N(0, 0.01^2)$, respectively.

[Table 2]

In this case, to provide references of posterior parameter distributions, DREAM$_{(ZS)}$ with three parallel chains and altogether 30,000 model evaluations were implemented. The likelihood for DREAM$_{(ZS)}$ was Gaussian likelihood with heteroscedastic measurement errors. The convergence was reached after about 20,000 model evaluations, and the last 10,000 samples were used to estimate the posterior distributions of parameters. For the IGPMC algorithm, the initial number of random parameter sets was 500, and 2,000 iterative steps were employed, i.e., the total number of model evaluation was 2,500. The trace plots of the 11 parameters generated by the IGPMC (i.e., all the candidate base points sequentially generated during the inverse GP construction) are shown in Figure 8. Clearly, the IGPMC algorithm quickly converges to the true values after the initial 500 random parameter sample drawings. After convergence, the IGPMC and DREAM$_{(ZS)}$ algorithms had rather similar probability distributions of the unknown parameters, which are shown in Figure 9. The true value, mean and standard deviation values obtained by the IGPMC algorithm for each parameter are also listed in Table 2.

[Figure 8]

[Figure 9]

In this case, for the IGPMC algorithm, although the number of unknown parameters is doubled (11 parameters in this case compared to 5 parameters in Case Study 2), the total number of model evaluations is similar to that in Case Study 2. While for MCMC algorithm, with more unknown parameters and increasing system complexity, many more model evaluations are needed to reach convergence.

To further improve the computational efficiency of IGPMC, parallel computation can be adopted. For the initial independent and random parameter sets drawn from prior distributions, the calculation of corresponding model outputs can be easily realized in a parallel mode. Meanwhile, in the subsequent iterative process, to better explore the propagation of measurement uncertainty to parameter uncertainty, one can generate multiple measurement realizations once, obtain multiple parameter estimates and calculate corresponding model outputs simultaneously in each iterative step, which is also in parallel. It should be noted that, MCMC algorithm can also adopt parallel computation to calculate the likelihood of parameter samples in the parallel chains.

We also used a limited number of model evaluations to demonstrate the efficiency of the IGPMC algorithm. The trace scatter plots for the parameter samples versus the true parameter values are shown in Figures 10 for this case. Here the minimum number of model evaluations was explored, where the single modality of the posterior distribution was assumed to be known. For this case, with 100 initial random parameter sets, less than 100 iterative steps were needed to obtain accurate parameter estimation. Under the same settings, this computational efficiency is hardly achieved by other parameter estimation methods. To make sure that the estimation results are reliable by IGPMC algorithm, more model evaluations is suggested in practical implementation.

[Figure 10]

## 3.4. Case Study 4 : A Conductivity Field Estimation Case

This case study tests the applicability of the IGPMC algorithm in inferring the spatially varying conductivity field, where the number of unknown parameters is much larger than the previous cases.

In this case, transient saturated flow was considered in a $800[L] \times 800[L]$ domain uniformly discretized into $41 \times 41$ grids. The upper and lower boundaries were no-flow, while the left and right boundaries were constant head boundaries with prescribed pressure heads of $202[L]$ and $198[L]$, respectively. The log conductivity field $Y(x)$ was modeled as a spatially correlated Gaussian random field with the following separable exponential correlation form,

$$C_Y(\mathbf{x}_1, \mathbf{x}_2) = C_Y(x_1, y_1; x_2, y_2) = \sigma_Y^2 \exp\left[-\frac{|x_1 - x_2|}{\lambda_x} - \frac{|y_1 - y_2|}{\lambda_y}\right], \qquad (8)$$

where $\sigma_Y^2 = 1$ is the variance and $\lambda_x = 320[L]$ and $\lambda_y = 320[L]$ are correlation lengths in the $x$ and $y$ directions, respectively. Then the Karhunen-Loève (KL) expansion was used to parameterize the log conductivity field. In this way, the number of unknown parameters was reduced from the total grid number (1681) to the truncated number of KL terms (40),

$$Y(\mathbf{x}) \approx \overline{Y}(\mathbf{x}) + \sum_{i=1}^{40} \xi_i \sqrt{\lambda_i} f_i(\mathbf{x}), \tag{9}$$

where $\overline{Y}(\mathbf{x})$ is the mean component, $\xi_i, i=1,2,\ldots,40$ are independent standard Gaussian random variables, $\lambda_i$ and $f_i(x)$ are eigenvalues and eigenfunctions of the covariance function described in Eq. (8). In this case study, the mean component $\overline{Y}(\mathbf{x})$ was zero. And the parameters to be estimated were the 40 standard Gaussian random variables, $\xi_i, i=1,2,\ldots,40$. The true conductivity field was generated with 100 KL terms, thus the model structural error was introduced. The transient hydraulic head measurements were generated with the reference field with additive Gaussian errors $\varepsilon \sim N(0, 0.01^2)$ at the 25 measurement locations every 0.6[T] up to 6[T], as shown in Figure 11(a).

[Figure 11]

Here, the prior distributions for the 40 variables in KL expansion were independent standard Gaussian distributions within the truncated bounds [-4, 4]. With the same hydraulic head measurements, the 40 Gaussian random variables in KL expansion were inferred with DREAM(ZS) and IGPMC algorithm, respectively. For DREAM(ZS), using Gaussian likelihood with homoscedastic measurement error, 3 parallel chains were used and totally 300,000 model evaluations were invoked, about 65,000 model simulations were needed for DREAM$_{(ZS)}$ to reach convergence. The sample of parameters that best fits the measurements in the 3 parallel chains (i.e., with the smallest RMSE value) was chosen to be the estimation of the log K field. For IGPMC, 1,000 initial random parameter samples were generated and 2,000 iteration steps were used, i.e., the total number of model evaluations was 3,000. Also, the parameter sample among all the base points that best fits the measurements was used to infer the log K field.

The RMSE values of parameter chains for DREAM$_{(ZS)}$ and sequential parameter base points for IGPMC are plotted in Figure 12 (a) and (b), respectively. The RMSE values for IGPMC decrease rapidly after the initial random parameter drawings. Although the RMSE values for IGPMC are with more variation, the smallest RMSE value (0.111) for the IGPMC algorithm is close (although not as good as) to that of DREAM$_{(ZS)}$ (0.101).

Figure 11 (b-d) shows the true log K field, the estimated fields given by DREAM$_{(ZS)}$ and IGPMC (with the smallest RMSE values), respectively. It clearly shows that both DREAM$_{(ZS)}$ and IGPMC can identify the main patterns of the true field. The pair wise comparisons between the true log K and the estimated values from DREAM$_{(ZS)}$ and IGPMC are shown in Figure 13. It shows that the estimated fields match the true field well.

[Figure 12]

[Figure 13]

It should be noted here that, the huge number of original model evaluations for MCMC algorithm could be greatly relieved though using surrogate systems. For example, using polynomial chaos expansion to construct a surrogate for the original system model, [*Laloy et al.*, 2013] developed a two-stage MCMC algorithm to efficiently explore the posterior of a high -dimensional groundwater model.

As shown in the above four case studies, IGPMC can obtain rather accurate parameter estimations for problems with bimodal distribution (Case Study 1), low (Case Study 1 and 2), moderate (Case Study 3) and high (Case Study 4) parameter dimensionality, respectively. While the computational cost for IGPMC is considerably low. As a matter of fact, the numbers of model evaluations are similar for different levels of dimensionality. For the four case studies tested in this paper, the numbers of model evaluations are all less than 3,000. Thus the new algorithm is especially efficient when the parameter dimensions are high.

## 4. Conclusions

An efficient parameter estimation and uncertainty assessment method, entitled the Inverse Gaussian Process Monte Carlo algorithm has been developed in this work. The IGPMC algorithm is conceptually simple and easy to implement. This method is based on iteratively constructing model output-input relationship through Gaussian Process regression and obtaining parameter estimate directly by using measurements as inputs to the inverse GP system. In this way, the uncertainty of measurements can be propagated to model parameters through the inverse GP system.

The efficiency and accuracy of the developed IGPMC algorithm in estimating hydrologic model parameters were tested in four numerical cases. In Case Study 1 with non-unique solutions, the bimodal probability distribution could be well identified by the new method. In Case Studies 2-4 with increasing dimension and complexity, the IGPMC algorithm could obtain accurate parameter estimations results with an affordable computational cost. Meanwhile, the time needed by the IGPMC algorithm could be further reduced through parallel computation.

# Acknowledgments.

Computer codes used are available upon request to the corresponding author.

We acknowledge Jasper Vrugt from University of California, Irvine for providing us with codes of DREAM$_{(ZS)}$.

Weixuan Li's work was supported by the Laboratory-Directed Research and Development program at Pacific Northwest National Laboratory (PNNL) through the Control of Complex Systems Initiative, and U.S. Department of Energy, Office of Science, Office of Advanced Scientific Computing Research, Applied Mathematics program as part of the Multifaceted Mathematics for Complex Energy Systems (M2ACS) project. PNNL is operated by Battelle Memorial Institute for the U. S. Department of Energy.

# References.

Andrieu, C., N. de Freitas, A. Doucet, and M. I. Jordan (2003), An introduction to MCMC for machine learning, *Mach. Learn.*, *50*(1-2), 5-43, doi: 10.1023/a:1020281327116.

Carrera, J., A. Alcolea, A. Medina, J. Hidalgo, and L. J. Slooten (2005), Inverse problem in hydrogeology, *Hydrogeol. J.*, *13*(1), 206-222, doi: 10.1007/s10040-004-0404-7.

Day, W. H. E., and H. Edelsbrunner (1984), Efficient algorithms for agglomerative hierarchical clustering methods, *J. Classif.*, *1*(1), 7-24, doi: 10.1007/bf01890115.

El Harrouni, K., D. Ouazar, G. A. Walters, and A. H. D. Cheng (1996), Groundwater optimization and parameter estimation by genetic algorithm and dual reciprocity boundary element method, *Eng. Anal. Bound. Elem.*, *18*(4), 287-296, doi: 10.1016/s0955-7997(96)00037-9.

Gupta, V. K., and S. Sorooshian (1985), The automatic calibration of conceptual catchment models using derivative - based optimization algorithms, *Water Resour. Res.*, *21*(4), 473-485, doi: 10.1029/WR021i004p00473.

Haario, H., M. Laine, A. Mira, and E. Saksman (2006), DRAM: Efficient adaptive MCMC, *Stat. Comput.*, *16*(4), 339-354, doi: 10.1007/s11222-006-9438-0.

Harbaugh, A. W., E. R. Banta, M. C. Hill, and M. G. McDonald (2000), *MODFLOW-2000, the US Geological Survey modular ground-water model: User guide to modularization concepts and the ground-water flow process*, US Geological Survey Reston, VA, USA.

Hartigan, J. A., and M. A. Wong (1979), Algorithm AS 136: A k-means clustering algorithm, *Appl. Stat.*, *28*(1), 100-108, doi: 10.2307/2346830

Jones, D. R., M. Schonlau, and W. J. Welch (1998), Efficient global optimization of expensive black-box functions, *J. Global Optim.*, *13*(4), 455-492, doi: 10.1023/A:1008306431147.

Kennedy, M. C., and A. O'Hagan (2001), Bayesian calibration of computer models, *J. R. Stat. Soc. Ser. B-Stat. Methodol.*, *63*(3), 425-450, doi: 10.1111/1467-9868.00294.

Laloy, E., B. Rogiers, J. A. Vrugt, D. Mallants, and D. Jacques (2013), Efficient posterior exploration of a high- dimensional groundwater model from two- stage Markov chain Monte Carlo simulation and polynomial chaos expansion, *Water Resour. Res.*, *49*(5), 2664-2682, doi: 10.1002/wrcr.20226.

Li, W. (2014), Inverse Modeling and Uncertainty Quantification of Nonlinear Flow in Porous Media Models, PhD thesis, University of Southern California, Los Angeles, California, USA.

Liu, X. Y., M. A. Cardiff, and P. K. Kitanidis (2010), Parameter estimation in nonlinear environmental problems, *Stoch. Environ. Res. Risk Assess.*, *24*(7), 1003-1022, doi: 10.1007/s00477-010-0395-y.

Ma, X., and N. Zabaras (2009), An efficient Bayesian inference approach to inverse problems based on an adaptive sparse grid collocation method, *Inverse Probl.*, *25*(3), 035013, doi: 10.1088/0266-5611/25/3/035013.

Matott, L. S., J. E. Babendreier, and S. T. Purucker (2009), Evaluating uncertainty in integrated environmental models: A review of concepts and tools, *Water Resour. Res.*, *45*(6), W06421, doi: 10.1029/2008wr007301.

Moore, R. J. (1985), The probability-distributed principle and runoff production at point and basin scales, *Hydrol. Sci. J.*, *30*(2), 273-297, doi: 10.1080/02626668509490989.

Oliver, D. S., and Y. Chen (2011), Recent progress on reservoir history matching: a review, *Computat. Geosci.*, *15*(1), 185-221, doi: 10.1007/s10596-010-9194-2.

Pan, L., and L. Wu (1998), A hybrid global optimization method for inverse estimation of hydraulic parameters: Annealing‐Simplex Method, *Water Resour. Res.*, *34*(9), 2261-2269, doi: 10.1029/98WR01672.

Rasmussen, C. E., and H. Nickisch (2010), Gaussian Processes for Machine Learning (GPML) Toolbox, *J. Mach. Learn. Res.*, *11*, 3011-3015.

Razavi, S., B. A. Tolson, and D. H. Burn (2012), Review of surrogate modeling in water resources, *Water Resour. Res.*, *48*(7), W07401, doi: 10.1029/2011wr011527.

Schoups, G., and J. A. Vrugt (2010), A formal likelihood function for parameter and predictive inference of hydrologic models with correlated, heteroscedastic, and non-Gaussian errors, *Water Resour. Res.*, *46*(10), W10531, doi: 10.1029/2009wr008933.

Smith, T. J., and L. A. Marshall (2008), Bayesian methods in hydrologic modeling: A study of recent advancements in Markov chain Monte Carlo techniques, *Water Resour. Res.*, *44*(12), W00B05, doi: 10.1029/2007wr006705.

Sun, A. Y., D. B. Wang, and X. L. Xu (2014), Monthly streamflow forecasting using Gaussian Process Regression, *J. Hydrol.*, *511*, 72-81, doi: 10.1016/j.jhydrol.2014.01.023.

Tartakovsky, D. M. (2013), Assessment and management of risk in subsurface hydrology: A review and perspective, *Adv. Water Resour.*, *51*, 247-260, doi: 10.1016/j.advwatres.2012.04.007.

Vrugt, J. A., H. V. Gupta, W. Bouten, and S. Sorooshian (2003), A Shuffled Complex Evolution Metropolis algorithm for optimization and uncertainty assessment of hydrologic model parameters, *Water Resour. Res.*, *39*(8), 1201, doi: 10.1029/2002wr001642.

Vrugt, J. A., C. J. F. ter Braak, H. V. Gupta, and B. A. Robinson (2009a), Equifinality of formal (DREAM) and informal (GLUE) Bayesian approaches in hydrologic modeling?, *Stoch. Environ. Res. Risk Assess.*, *23*(7), 1011-1026, doi: 10.1007/s00477-008-0274-y.

Vrugt, J. A., C. J. F. ter Braak, M. P. Clark, J. M. Hyman, and B. A. Robinson (2008), Treatment of input uncertainty in hydrologic modeling: Doing hydrology backward with Markov chain Monte Carlo simulation, *Water Resour. Res.*, *44*, W00B09, doi: 10.1029/2007wr006720.

Vrugt, J. A., C. J. F. ter Braak, C. G. H. Diks, B. A. Robinson, J. M. Hyman, and D. Higdon (2009b), Accelerating Markov Chain Monte Carlo Simulation by Differential Evolution with Self-Adaptive Randomized Subspace Sampling, *Int. J. Nonlin. Sci. Num.*, *10*(3), 273-290, doi: 10.1515/IJNSNS.2009.10.3.273.

Woodbury, A. D., and T. J. Ulrych (2000), A full-Bayesian approach to the groundwater inverse problem for steady state flow, *Water Resour. Res.*, *36*(8), 2081-2093, doi: 10.1029/2000wr900086.

Zeng, L., L. Shi, D. Zhang, and L. Wu (2012), A sparse grid based Bayesian method for contaminant source identification, *Adv. Water Resour.*, *37*, 1-9, doi: 10.1016/j.advwatres.2011.09.011.

Zhang, G., D. Lu, M. Ye, M. Gunzburger, and C. Webster (2013), An adaptive sparse-grid high-order stochastic collocation method for Bayesian inference in groundwater reactive transport modeling, *Water Resour. Res.*, *49*(10), 6871-6892, doi: 10.1002/wrcr.20467.

Zhang, J., L. Zeng, C. Chen, D. Chen, and L. Wu (2015), Efficient Bayesian experimental design for contaminant source identification, *Water Resour. Res.*, *51*(1), 576-598, doi: 10.1002/2014WR015740.

Zheng, C., and P. P. Wang (1999), MT3DMS: a modular three-dimensional multispecies transport model for simulation of advection, dispersion, and chemical reactions of contaminants in groundwater systems; documentation and user's guide*Rep.*, DTIC Document.

Zhou, H., J. Jaime Gomez-Hernandez, and L. Li (2014), Inverse methods in hydrogeology: Evolution and recent trends, *Adv. Water Resour.*, *63*, 22-37, doi: 10.1016/j.advwatres.2013.10.014.

## Table Captions:

**Table 1.** Prior range, true value, mean (Mean) and standard deviation (SD) values obtained by the IGPMC algorithm for each HYMOD model parameter.

**Table 2.** Prior range, true value, mean (Mean) and standard deviation (SD) values obtained by the IGPMC algorithm for each unknown parameter in Case Study 3.

# Figure Captions:

**Figure 1.** Schematic diagram for the Inverse Gaussian Process Monte Carlo algorithm.

**Figure 2.** Bimodal probability distribution of parameter in Case Study 1 inferred from DREAM$_{(ZS)}$ (represented by blue line) and the IGPMC algorithm (represented by red dashed line). The true value of parameter is represented by vertical black line.

**Figure 3.** Trace plots of (a) $C_{max}$, (b) $b_{exp}$, (c) $\alpha*$, (d) $R_s$, (e) $R_q$ for all the candidate base points generated through the implementation of the IGPMC algorithm.

**Figure 4.** Probability distributions of the HYMOD model parameters inferred with DREAM$_{(ZS)}$ (represented by blue lines) and the IGPMC algorithm (represented by red dashed lines). The true values of parameters are represented by vertical black lines.

**Figure 5.** Probability distributions of the HYMOD model parameters inferred with DREAM(ZS) (represented by blue lines) and the IGPMC algorithm (represented by red dashed lines) without knowing the measurement error statistics a priori. The true values of parameters are represented by vertical black lines.

**Figure 6.** Comparison between actual measurement errors and estimated measurement errors for Case Study 2.

**Figure 7.** Flow domain for Case Study 3.

**Figure 8** Trace plots of (a, b) source location parameters, (c-h) source strength parameters and (i-k) log conductivity parameters in Case Study 3 for all the candidate base points generated through the implementation of the IGPMC algorithm.

**Figure 9.** Probability distributions of contaminant transport model parameters inferred with DREAM$_{(ZS)}$ (represented by blue lines) and the IGPMC algorithm (represented by red dashed lines). The true values are represented by vertical black lines.

**Figure 10.** Trace plots of (a, b) source location parameters, (c-h) source strength

parameters and (i-k) log conductivity parameters in Case Study 3 for all the candidate base points generated through the implementation of the IGPMC algorithm with 200 model evaluations in total.

**Figure 11.** (a) The flow domain and measurement locations for the pressure head (25 filled dots); (b) True log K field; (c) Log K field estimated with DREAM$_{(ZS)}$; (d) Log K field estimated with the IGPMC algorithm.

**Figure 12.** The RMSE values between model outputs at all the parameter samples and measurement for (a) DREAM$_{(ZS)}$ and (b) the IGPMC algorithm.

**Figure 13.** Pair wise comparison of the true log K and estimated log K obtained through (a) DREAM$_{(ZS)}$ and (b) the IGPMC algorithm.

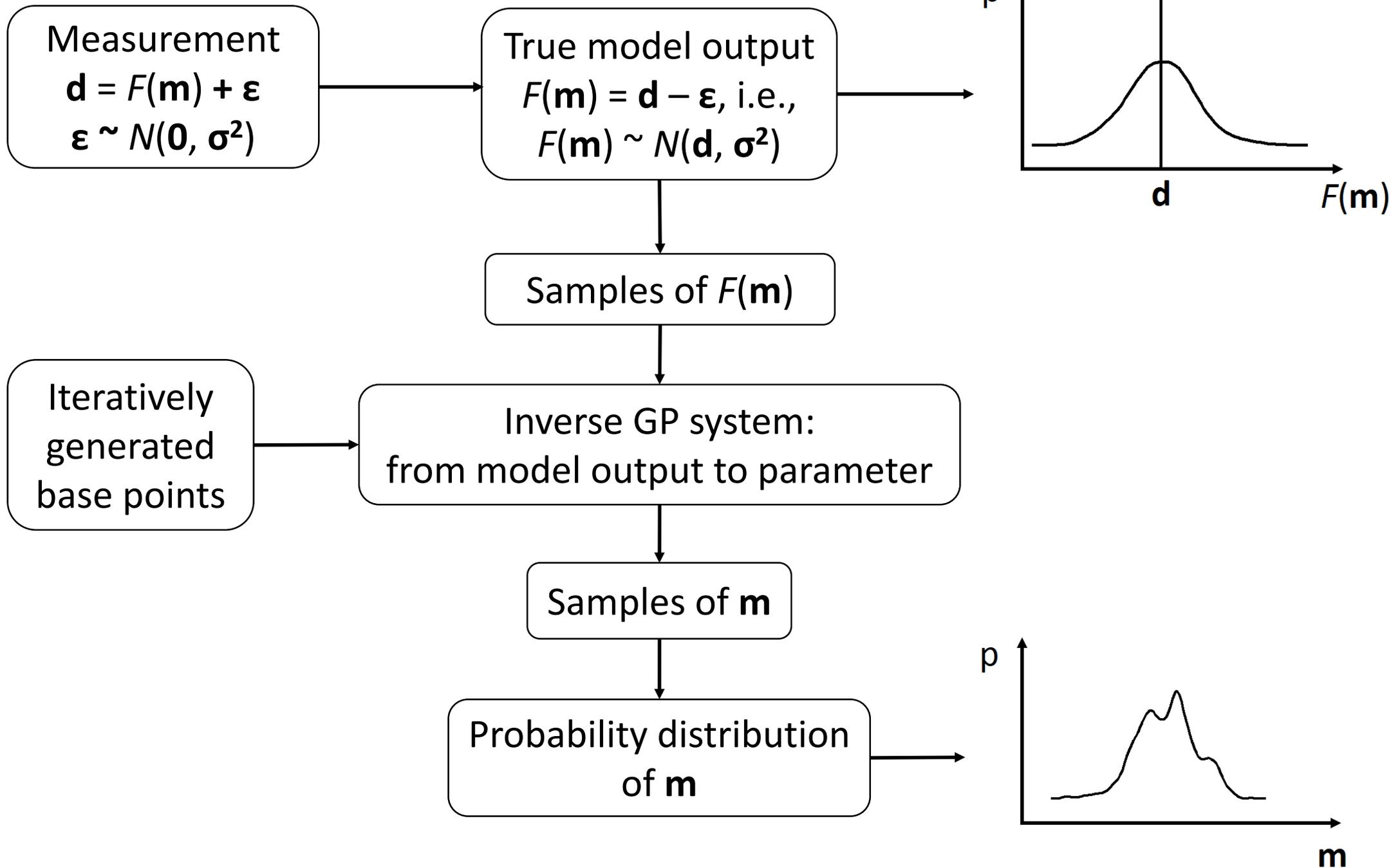

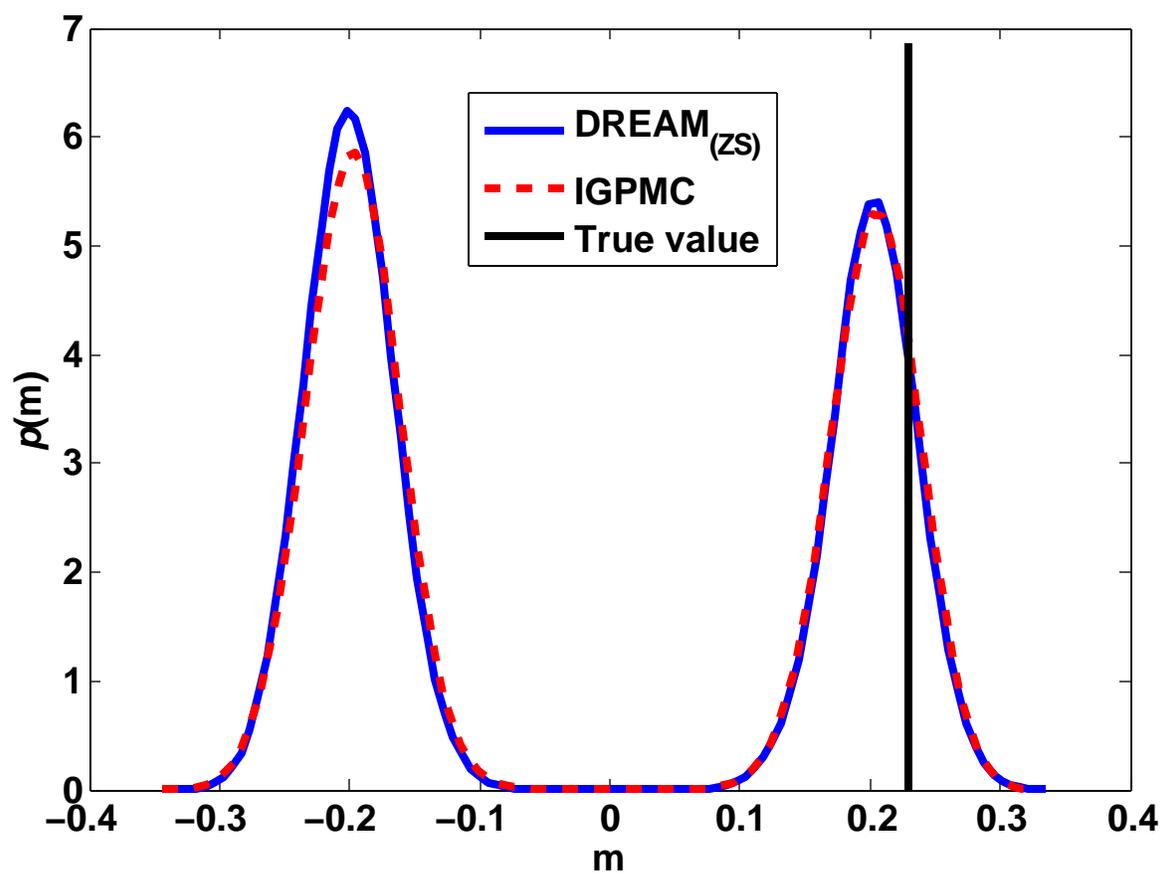

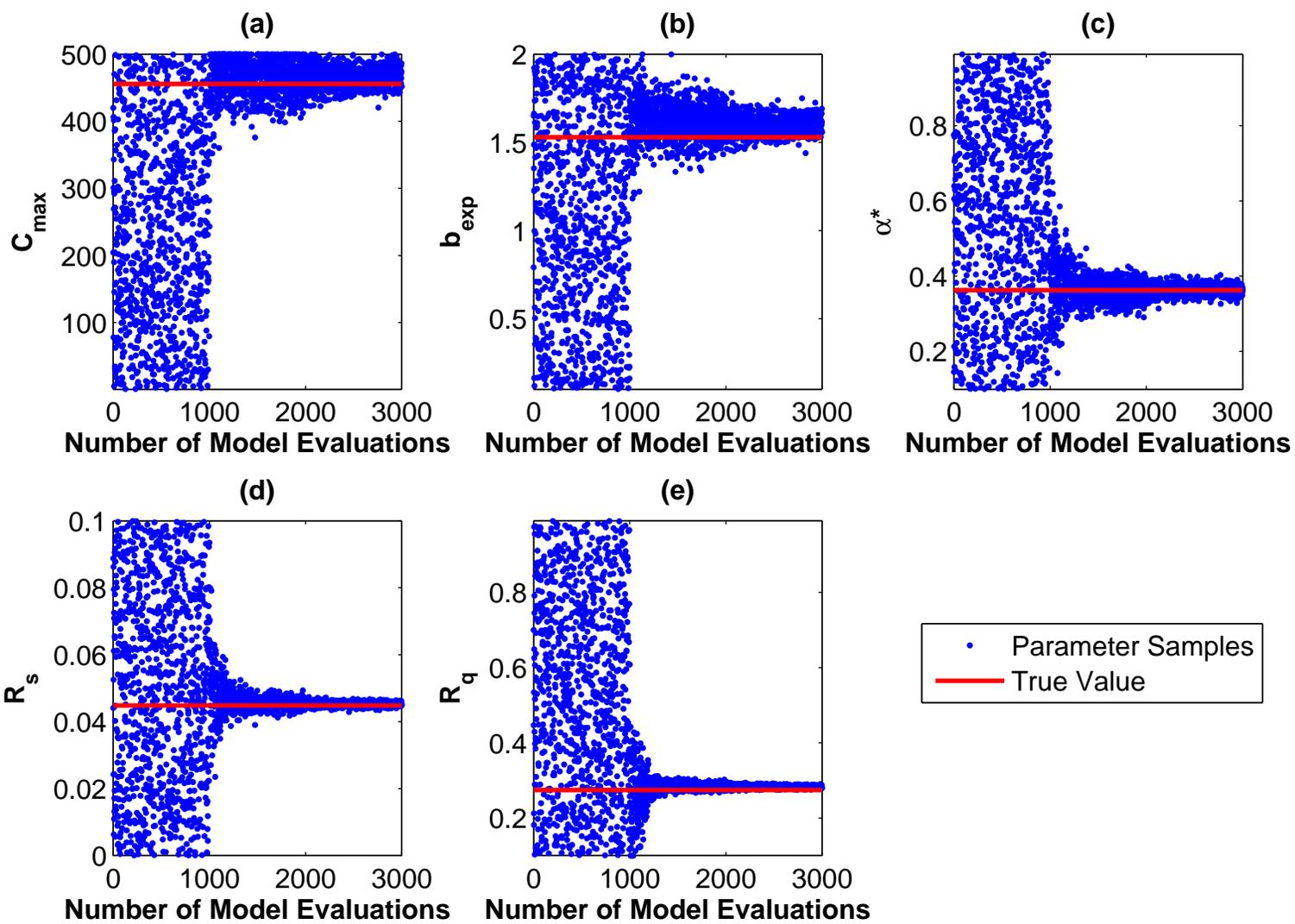

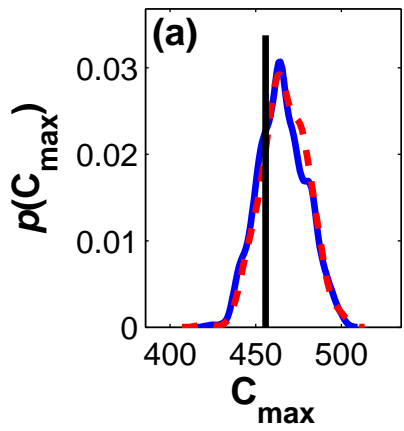
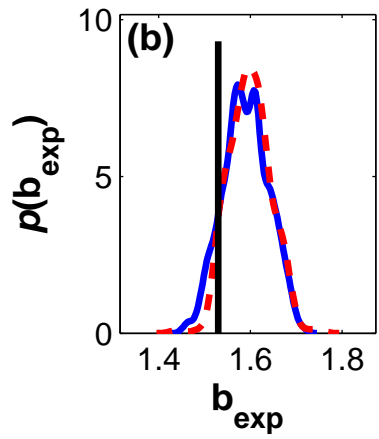
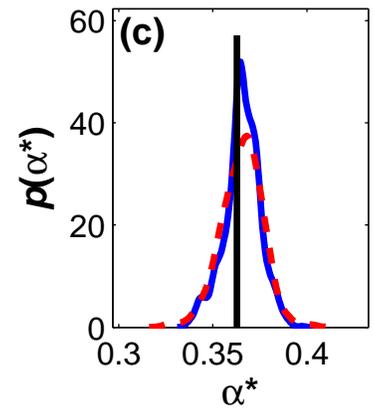
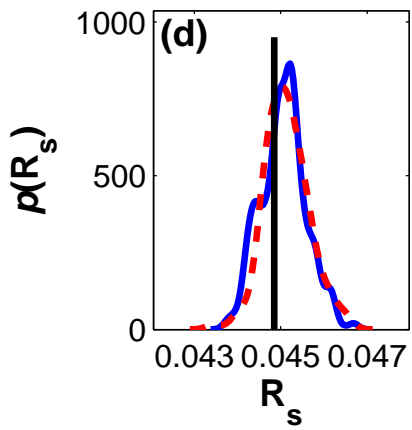
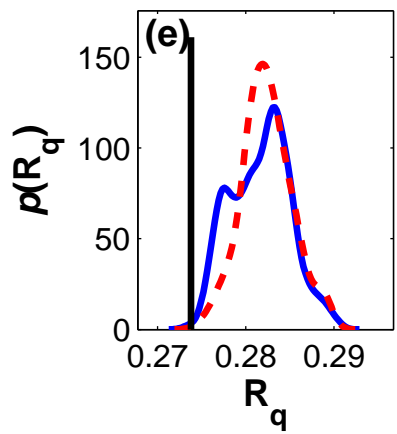
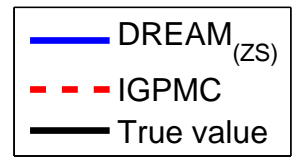

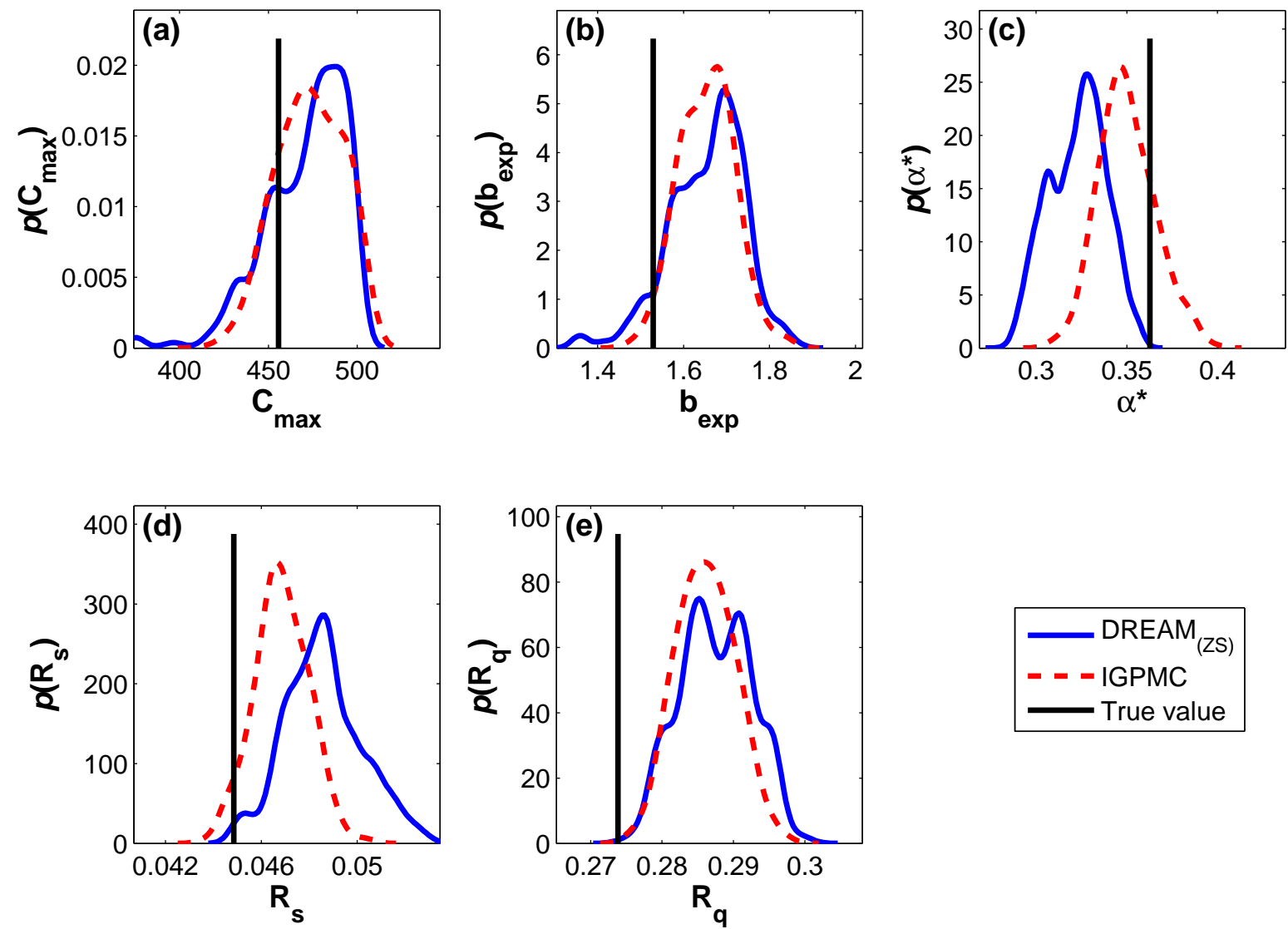

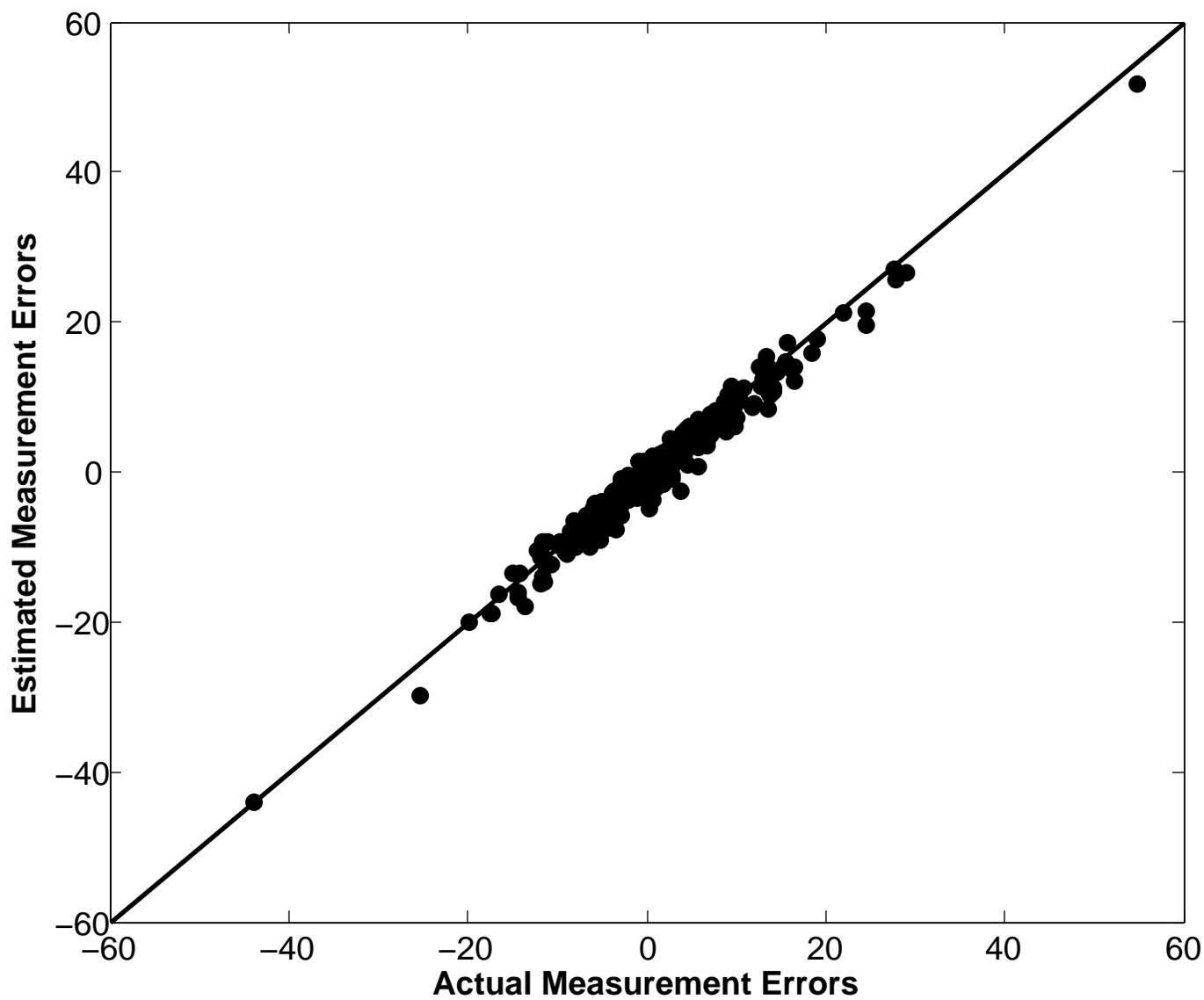

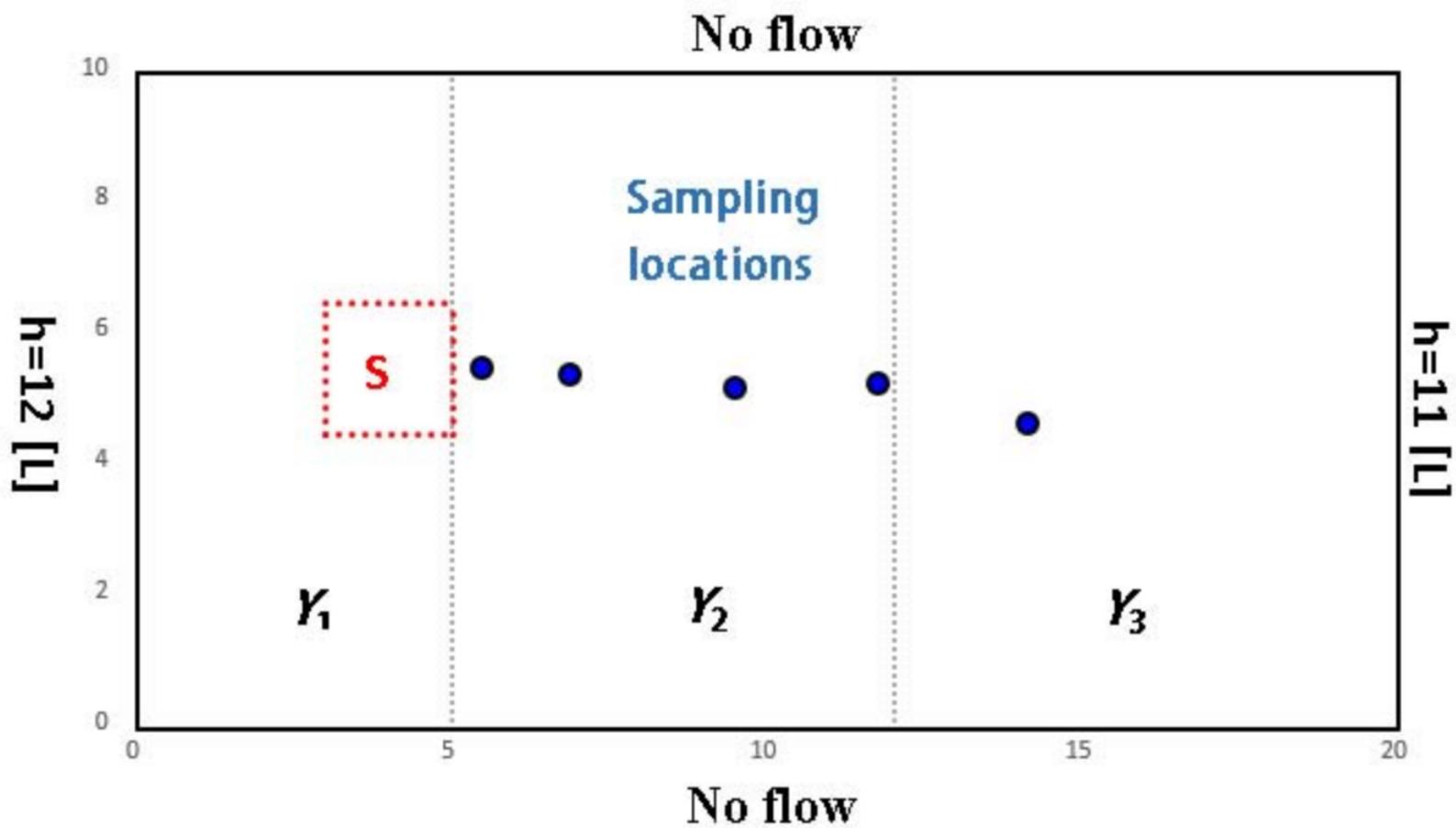

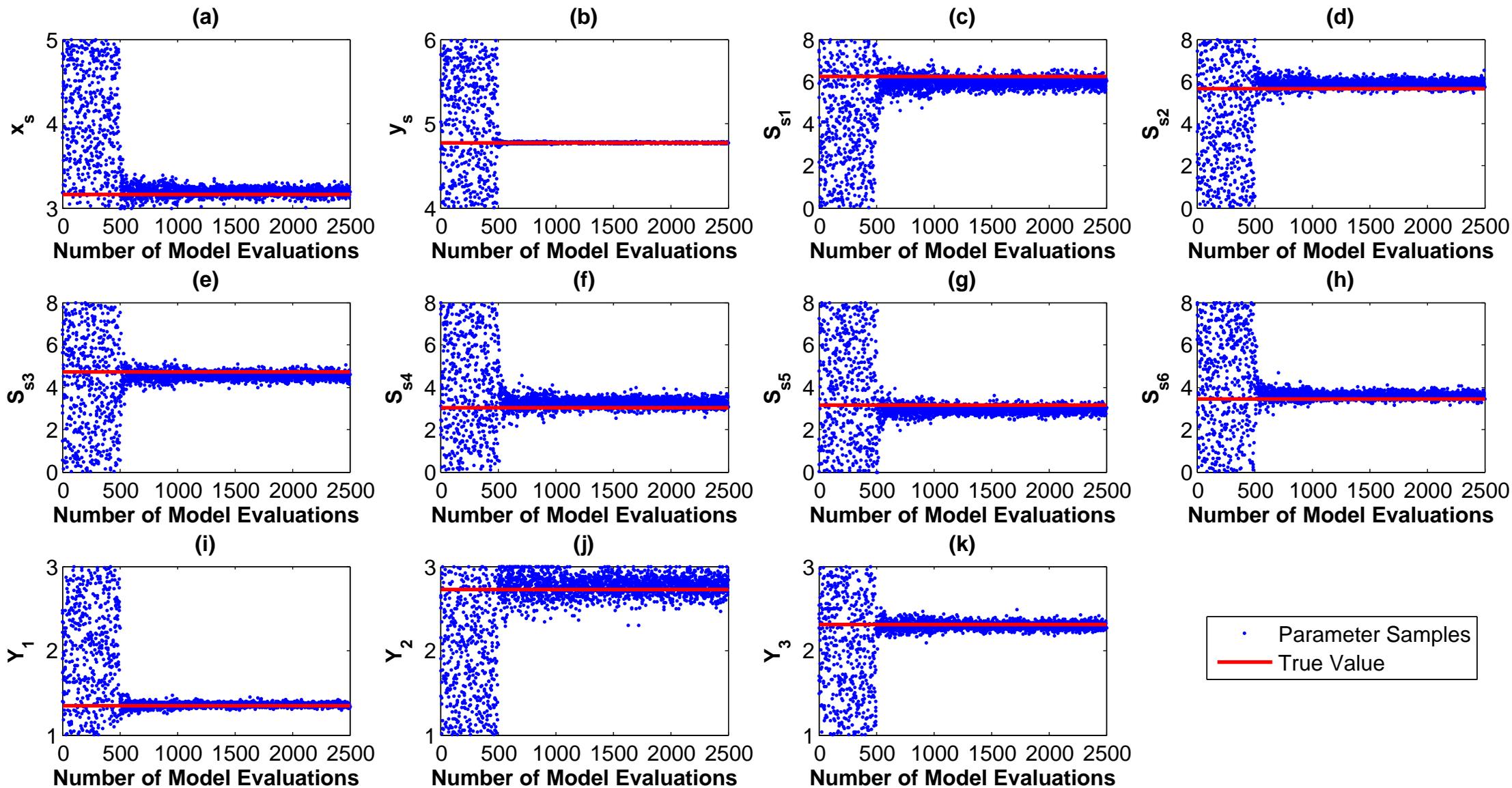

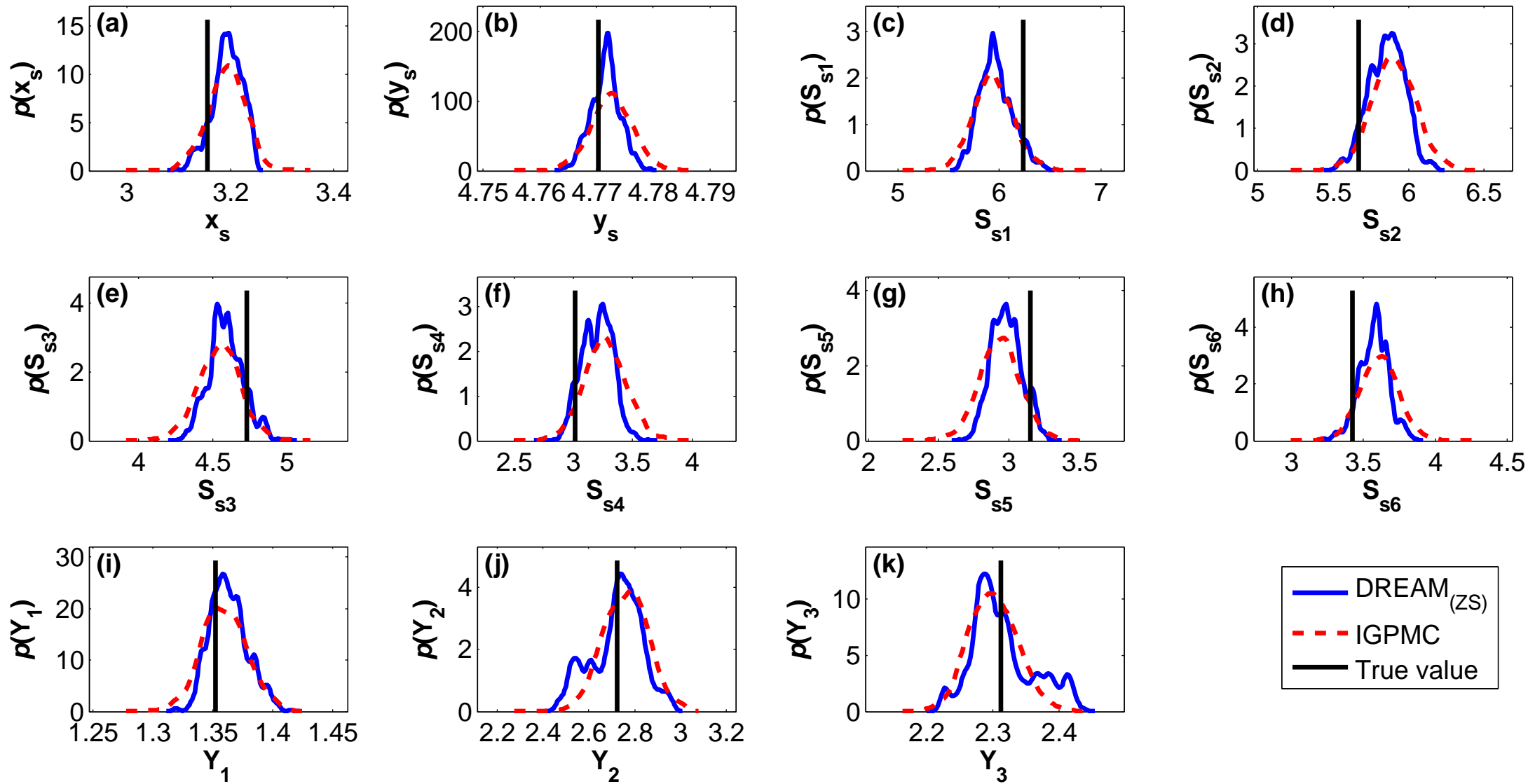

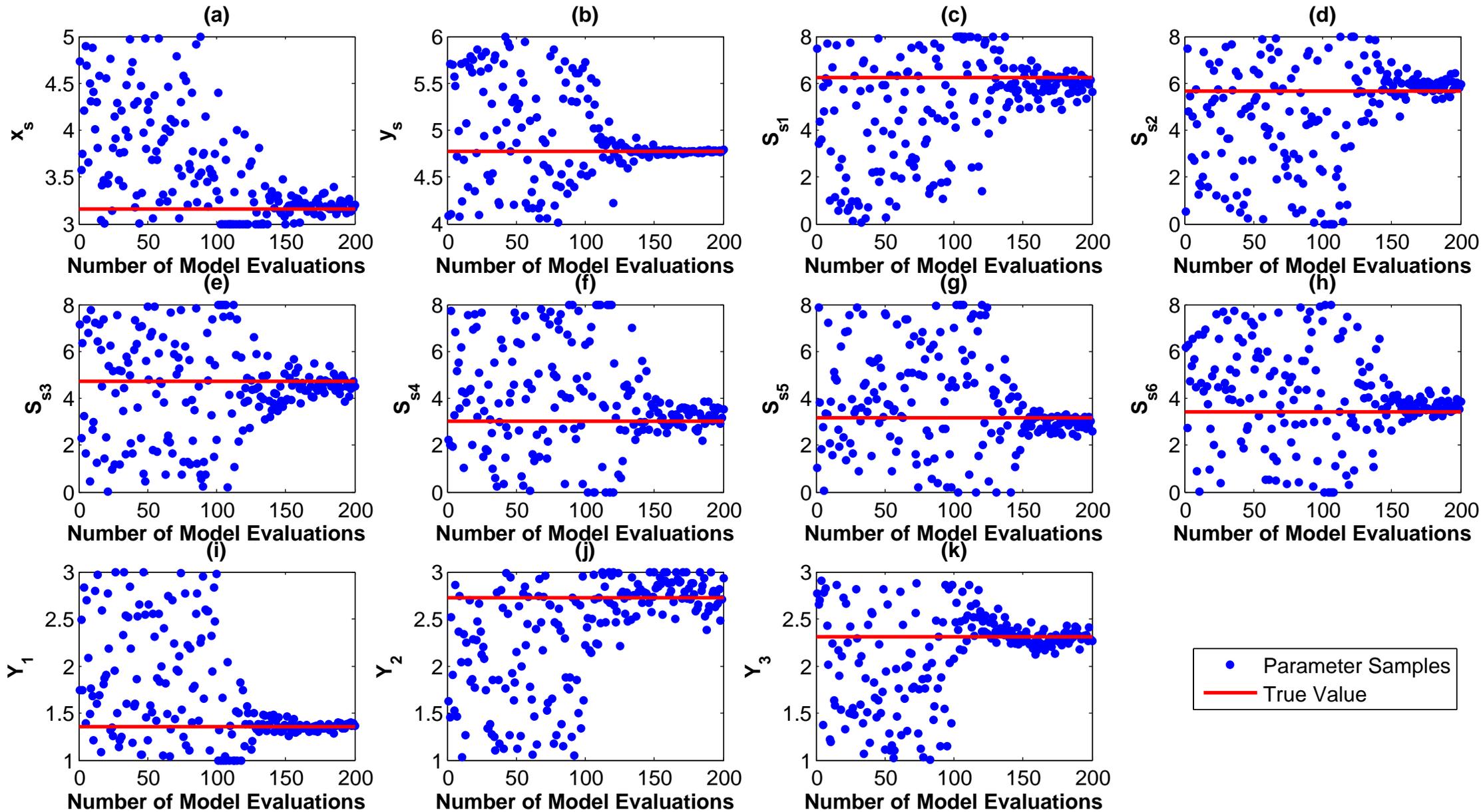

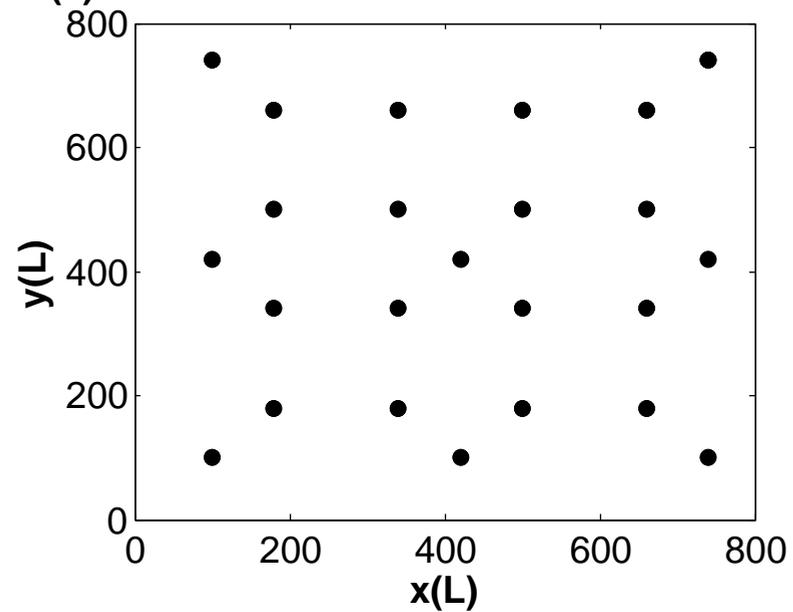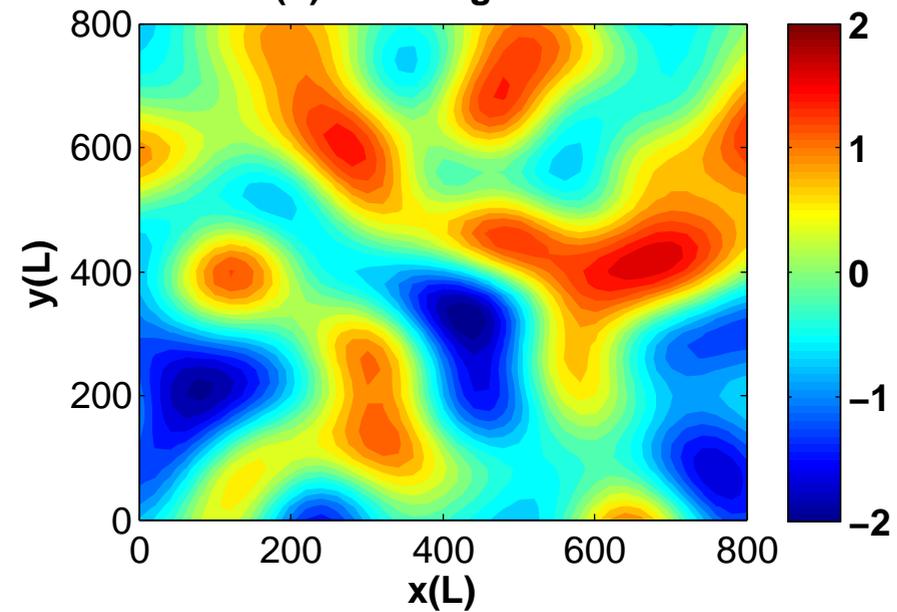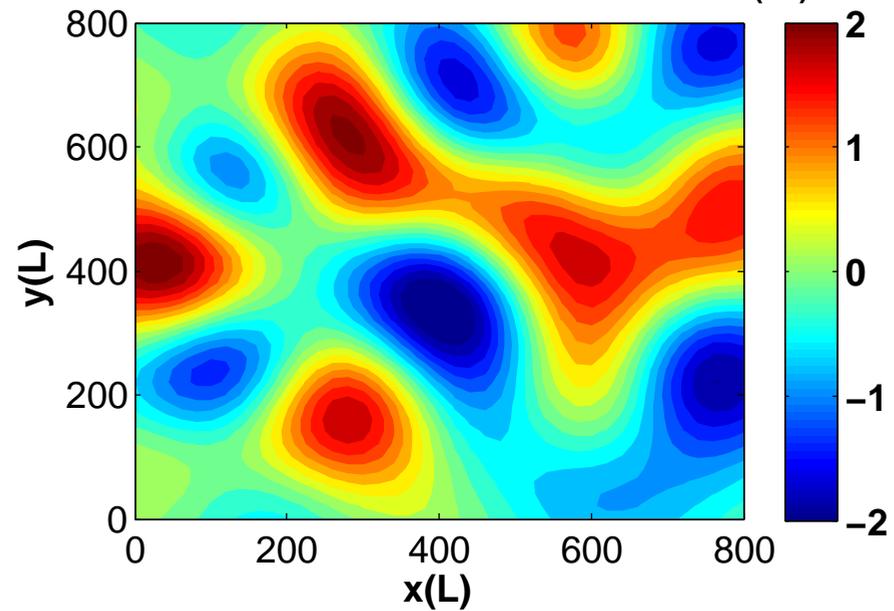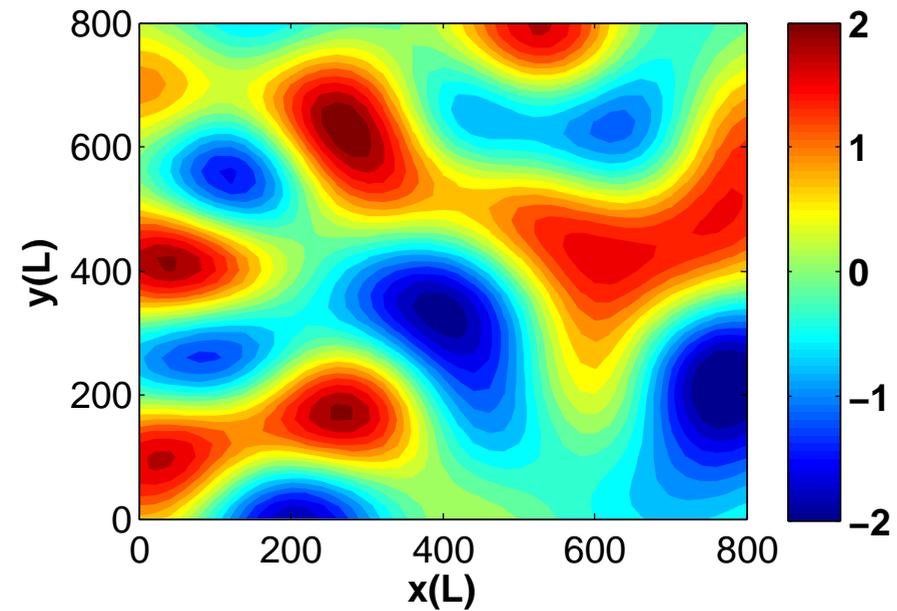

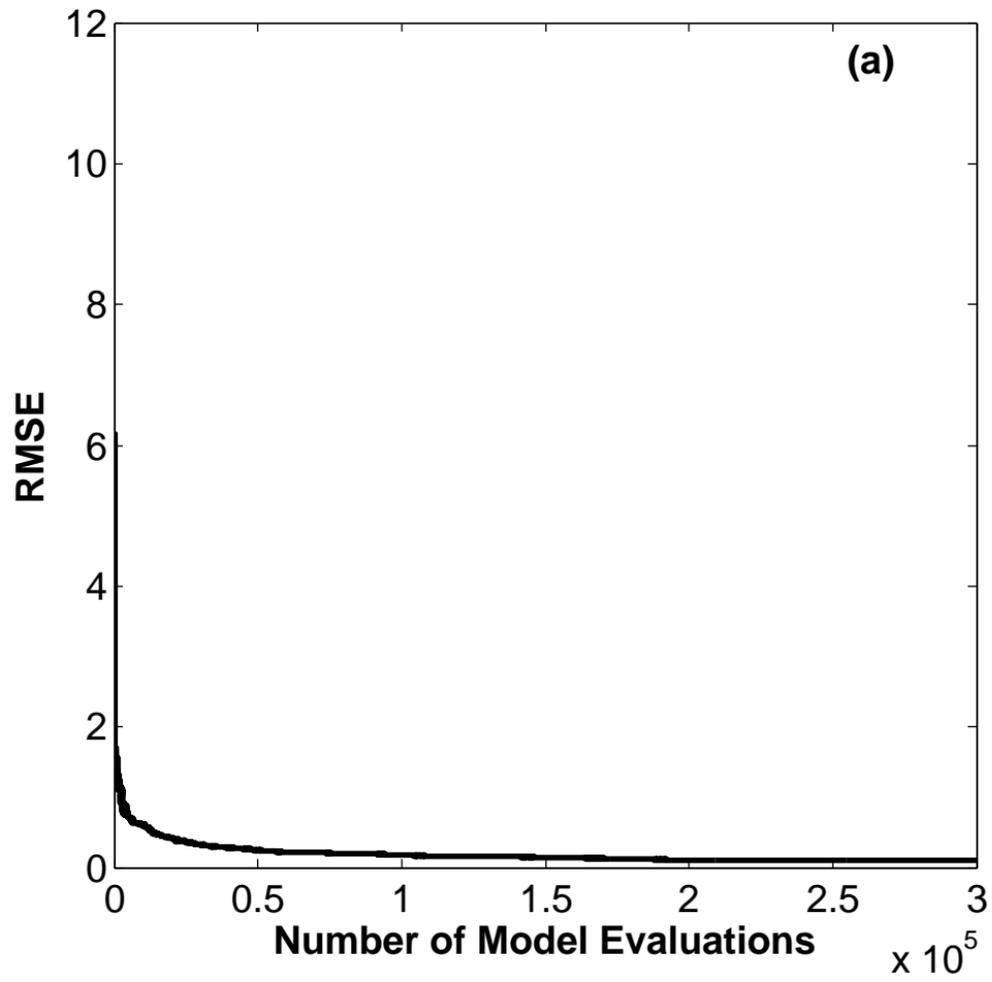 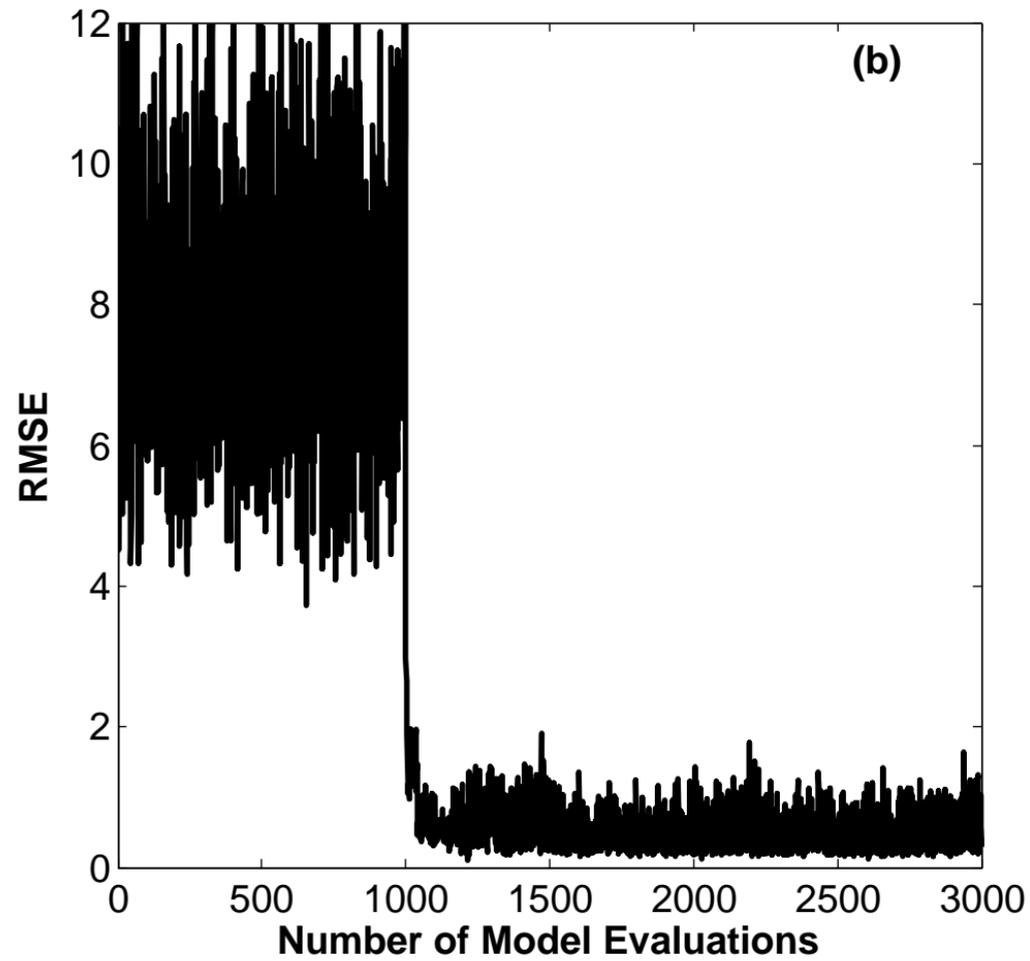

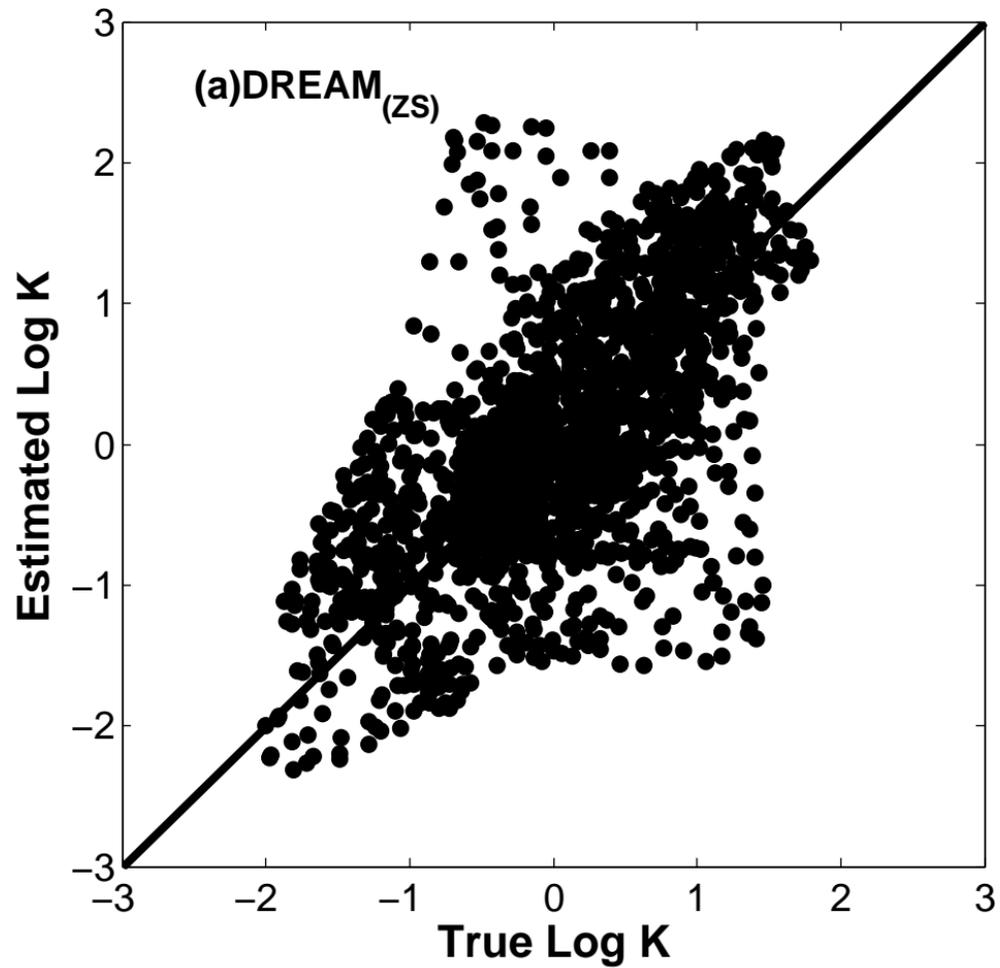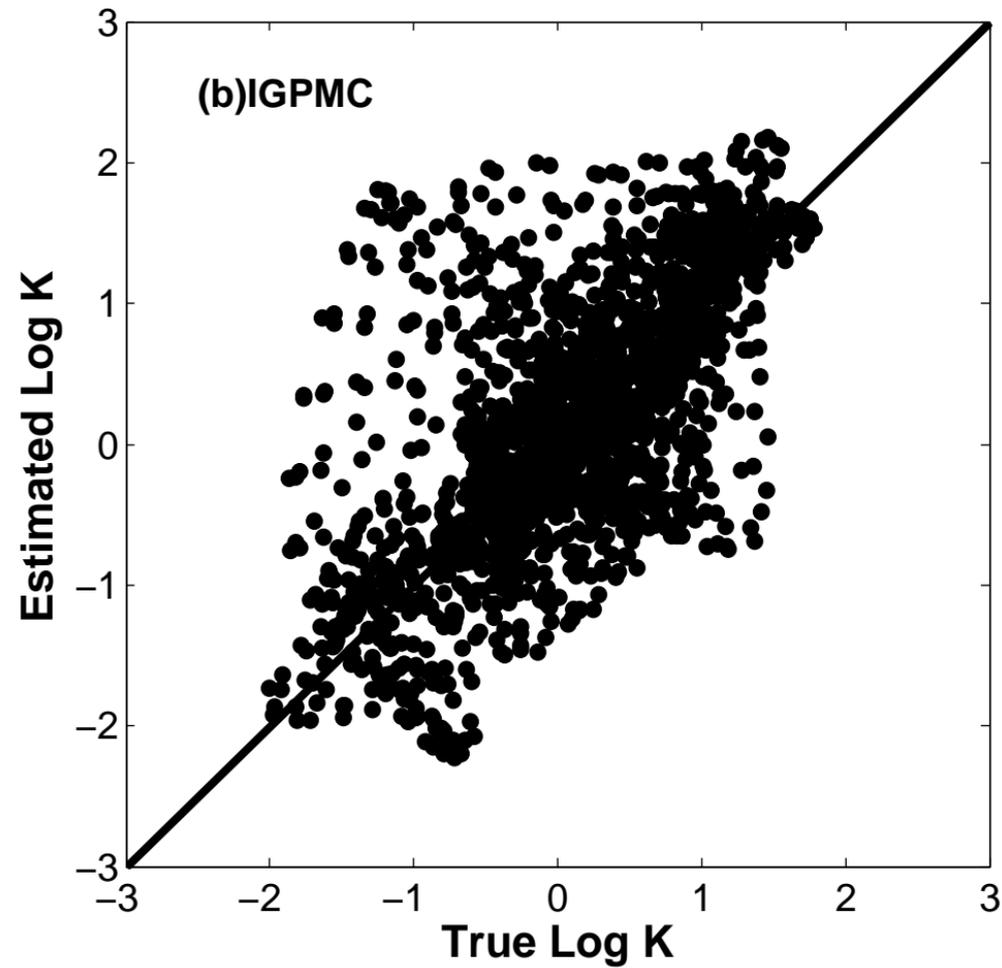